\def\mod{\mathop{\rm mod} \nolimits}
\baselineskip=14pt
\parskip=10pt
\def\Tilde{\char126\relax}
\def\halmos{\hbox{\vrule height0.15cm width0.01cm\vbox{\hrule height
 0.01cm width0.2cm \vskip0.15cm \hrule height 0.01cm width0.2cm}\vrule
 height0.15cm width 0.01cm}}
\font\eightrm=cmr8  
\font\eighttt=cmtt8
\magnification=\magstephalf

\parindent=0pt
\overfullrule=0in
\bf
\centerline
{An Explicit Formula for the Number of Solutions of $X^2=0$ in Triangular
Matrices Over a Finite Field}
\rm
\bigskip
\centerline{ {\it Shalosh B. EKHAD}{$^1$} and
{\it Doron ZEILBERGER}\footnote{$^1$}
{\eightrm  \raggedright
Department of Mathematics, Temple University,
Philadelphia, PA 19122, USA. 
{\eighttt [ekhad,zeilberg]@math.temple.edu,
http://www.math.temple.edu/\Tilde [ekhad, zeilberg] \quad
 ftp://ftp.math.temple.edu/pub/[ekhad,zeilberg]} \quad .
\break
Supported in part by the NSF. Nov. 28, 1995.
} 
}
{\bf Abstract:} We prove an explicit formula for the number of $n \times n$
upper triangular matrices, over $GF(q)$, 
whose square is the zero matrix. This formula was recently conjectured
by Sasha Kirillov and Anna Melnikov[KM].
 
{\bf Theorem}: The number of $n \times n$ upper-triangular matrices
over $GF(q)$ (the finite field with $q$ elements), whose square
is the zero matrix, is given by the polynomial $C_n(q)$, where,
$$
C_{2n}(q)=\sum_{j} \left [{{2n} \choose {n-3j}}-  
             {{2n} \choose {n-3j-1}}  \right ] \cdot q^{n^2-3j^2-j} \quad ,
$$
$$
C_{2n+1}(q)=\sum_{j} \left [{{2n+1} \choose {n-3j}}-  
             {{2n+1} \choose {n-3j-1}}  \right ] \cdot q^{n^2+n-3j^2-2j}
 \quad .
$$
 
{\bf Proof:} In [K] it was shown that the quantity of interest is
given by the polynomial $A_n(q)=\sum_{r \geq 0} A_n^{r}(q)$, where
the polynomials $A_n^r(q)$ are defined recursively by:
$$
A_{n+1}^{r+1}(q)=q^{r+1} \cdot A_{n}^{r+1}(q) +
(q^{n-r}-q^r) \cdot  A_{n}^{r}(q) \quad ; \quad A_{n+1}^{0}(q)=1 \quad.
\eqno(Sasha)
$$
 
For any Laurent formal power series $P(w)$, let
$CT_w P(w)$ denote the coefficient of $w^0$. 
Recall that the {\it $q$-binomial coefficients} are defined by
$$
{{m} \choose {n}}_q :=
{{(1-q^m)(1-q^{m-1}) \cdots (1-q^{m-n+1})} 
\over {(1-q)(1-q^2) \cdots (1-q^n)}} \quad ,
\eqno(Carl)
$$
whenever $0 \leq n \leq m$, and $0$ otherwise.
 
The following lemma gives an explicit expression for
$A_n^r(q)$.
 
{\bf Lemma 1:}
$$
A_n^r(q)=CT_w \left [
{{(1-w)(1+w)^n q^{r(n-r)}} \over {w^r}} 
\sum_{i=0}^{\infty} (-1)^i q^{-(i+1)i/2 -i(n-2r)}
{{i+n-2r} \choose {i}}_q  w^i \right ] \quad .
\eqno(Anna)
$$
 
{\bf Proof}: Call the right side of Eq. $(Anna)$, $S_n^r(q)$.
Since $S_{n+1}^{0}(q)=1$, the lemma would follow by induction if
we could show that
$$
S_{n+1}^{r+1}(q) \,- \,q^{r+1} \cdot S_{n}^{r+1}(q) \,-\,
(q^{n-r}-q^r) \cdot S_{n}^{r}(q)\,=0 \quad .
\eqno(Sasha')
$$
Using the linearity of $CT_w$, manipulating the series, using
the definition $(Carl)$ of the $q-$binomial coefficients, and
simplifying, brings the left side of $(Sasha')$ to be
$CT_w{\Phi_n^r(q,w)}$, where $\Phi_n^r$ is zero except when
$n$ is odd and $r=(n-1)/2$, in which case it is a monomial in $q$ times
${{(1-w)(1+w)^n} \over {w^{r+1}}} $, and applying $CT_w$ kills it
all the same, thanks to the symmetry of the Chu-Pascal triangle. \halmos
 
Summing the expression proved for $A_n^r(q)$, yields that
$$
A_n(q)=CT_w \left [(1-w)(1+w)^n \cdot
\sum_{r=0}^{\infty} \sum_{i=0}^{r}
 (-1)^i q^{r(n-r)-(i+1)i/2 -i(n-2r)}
{{i+n-2r} \choose {n-2r}}_q  w^{i-r} \right ] \quad .
$$
Letting $l=r-i$, and changing the order of summation, yields
$$
A_n(q)=CT_w \left [(1-w)(1+w)^n \cdot
\sum_{l=0}^{\lfloor n/2 \rfloor } w^{-l} \cdot q^{ln-l^2} \, 
\sum_{i=0}^{\lfloor (n-2l)/2 \rfloor} (-1)^i q^{i(i-1)/2}
{{n-2l-i} \choose {i}}_q \right ] \quad .
\eqno(SumAnna)
$$
Luckily, the inner sum evaluates nicely thanks to:
 
{\bf Lemma 2:}
$$
\sum_{i=0}^{\lfloor m/2 \rfloor}
 (-1)^i q^{i(i-1)/2}
{{m-i} \choose {i}}_q =(-1)^{\lfloor m/3 \rfloor} q^{m(m-1)/6} \cdot
\chi (m \not\equiv 2 \mod 3 ) \quad .
$$
{\bf Proof:} While this is unlikely to be 
new\footnote{$^2$}{\eightrm For a `classical' proof,
see Christian Krattenthaler's message, 
at {\eighttt ftp://ftp.math.temple.edu/pub/ekhad/sasha.}}, it is also 
irrelevant whether or not it is new, since
this is {\it now} routine, thanks to the
package {\tt qEKHAD}, accompanying [PWZ]. 
Let's call the left side divided by $q^{m(m-1)/6}$, $Z(m)$.
Then we have to prove that $Z_0(m):=Z(3m)$ equals $(-1)^m$,
$Z_1(m):=Z(3m+1)$ equals $(-1)^m$, and
$Z_2(m):=Z(3m+2)$ equals $0$. It is directly verified that
these are true for $m=0,1$, and the general result follows from
the second order recurrences produced by {\tt qEKHAD}.
The input files
{\tt inZ0}, {\tt inZ1}, {\tt inZ2} as well as the corresponding
output files, {\tt outZ0}, {\tt outZ1}, {\tt outZ2}
can be obtained by anonymous
{\tt ftp} to {\tt ftp.math.temple.edu}, directory {\tt pub/ekhad/sasha}.
The package {\tt qEKHAD} can be downloaded from
{\tt http://www.math.temple.edu/\Tilde zeilberg}. \halmos
 
To complete the proof of the theorem, we use lemma $2$ to 
evaluate the inner sum of $(SumAnna)$, then to get
$A_{2n}(q)$, we replace $n$ by $2n$, and then
replace $l$ by $l+n$, and finally use the binomial theorem.
Similarly for $A_{2n+1}(q)$. \halmos     
 
{\bf References}
 
[K] A.A. Kirillov, {\it On the number of solutions to the
equation $X^2=0$ in triangular matrices over a finite field},
Funct. Anal. and Appl. {\bf 29} (1995), no. 1.
 
[KM] A.A. Kirillov and A. Melnikov, {\it On a remarkable sequence
of polynomials}, preprint.
 
[PWZ] M. Petkovsek, H.S. Wilf, and D. Zeilberger, {\it ``A=B''},
A.K.Peters, 1996. 
 
\bye